\documentclass[11pt,reqno,letterpaper]{amsart}
\usepackage{amsmath,amsthm,amsfonts,amssymb,bbm}

\usepackage{bookmark}
\usepackage{hyperref}
\hypersetup{pdfstartview={FitH}}

\newtheorem{theorem}{Theorem}

\numberwithin{equation}{section}

\begin{document}

\title[Asymptotically sharp discrete nonlinear Hausdorff--Young inequalities]{Asymptotically sharp discrete nonlinear Hausdorff--Young inequalities for the $\textup{SU}(1,1)$-valued Fourier products}

\author[V. Kova\v{c}]{Vjekoslav Kova\v{c}}
\address{Vjekoslav Kova\v{c}, Department of Mathematics, Faculty of Science, University of Zagreb, Bijeni\v{c}ka cesta 30, 10000 Zagreb, Croatia}
\email{vjekovac@math.hr}

\author[D. Oliveira e Silva]{Diogo Oliveira e Silva}
\address{Diogo Oliveira e Silva, Departamento de Matemática, Instituto Superior Técnico,  Av. Rovisco Pais, 1049-001 Lisboa, Portugal}
\email{diogo.oliveira.e.silva@tecnico.ulisboa.pt}

\author[J. Rup\v{c}i\'{c}]{Jelena Rup\v{c}i\'{c}}
\address{Jelena Rup\v{c}i\'{c}, Faculty of Transport and Traffic Sciences, University of Zagreb, Vukeli\'{c}eva 4, 10000 Zagreb, Croatia}
\email{jrupcic@fpz.hr}




\begin{abstract}
We work in a discrete model of the nonlinear Fourier transform (following the terminology of Tao and Thiele), which appears in the study of orthogonal polynomials on the unit circle. The corresponding nonlinear variant of the Hausdorff--Young inequality can be deduced by adapting the ideas of Christ and Kiselev to the present discrete setting. However, the behavior of sharp constants remains largely unresolved. In this short note we give two results on these constants, after restricting our attention to {either} sufficiently small sequences or to sequences that are far from being the extremizers of the linear Hausdorff--Young inequality.
\end{abstract}

\maketitle


\section{Introduction and statement of the results}
Let $F=(F_n)_{n\in\mathbb{Z}}$ be a sequence of complex numbers in the unit disk $\{z\in\mathbb{C}:|z|<1\}$, with only finitely many $F_n$ being nonzero.
We treat them as certain coefficients and set
\[ A_n := \frac{1}{(1-|F_n|^2)^{1/2}},\quad B_n := \frac{F_n}{(1-|F_n|^2)^{1/2}}, \]
so that for each $n$ we indeed have $A_n>0$, $B_n\in\mathbb{C}$, $A_n^2-|B_n|^2=1$.
Alternatively, we could have started with $A_n$, $B_n$ satisfying these constraints and set $F_n=B_n/A_n$. For any $t\in\mathbb{T}=\mathbb{R}/\mathbb{Z}\equiv[0,1\rangle$ we define the \emph{$\textup{SU}(1,1)$-valued trigonometric (or Fourier) product with coefficients $F$} as the matrix product
\begin{align}
\begin{bmatrix} a(t) & b(t) \\ \overline{b(t)} & \overline{a(t)} \end{bmatrix}
& := \prod_{n=-\infty}^{\infty} \begin{bmatrix} A_n & B_n e^{2\pi i n t} \\ \overline{B_n} e^{-2\pi i n t} & A_n \end{bmatrix} \nonumber \\
& = \Big(\prod_{n\in\mathbb{Z}}(1-|F_n|^2)^{-1/2}\Big) \prod_{n=-\infty}^{\infty} \begin{bmatrix} 1 & F_n e^{2\pi i n t} \\ \overline{F_n} e^{-2\pi i n t} & 1 \end{bmatrix}. \label{eq:defforNFT}
\end{align}
Its name emphasizes the matrix group
\[ \textup{SU}(1,1) := \bigg\{ \begin{bmatrix} \alpha & \beta \\ \overline{\beta} & \overline{\alpha} \end{bmatrix} : \alpha,\beta\in\mathbb{C},\ |\alpha|^2 - |\beta|^2 = 1 \bigg\} \subset \textup{SL}(2,\mathbb{C}) \]
and the above product takes values in that group for any fixed $t$. In particular, $|a(t)|\geq 1$ for each $t\in\mathbb{T}$.
Matrix multiplication is performed from left to right as $n$ increases and the order is important.
Only finitely many matrices in the product are different from the identity matrix by our assumption on $F$, so the multiplication process is effectively finite and we do not need to discuss convergence.

This setup was suggested by Tao and Thiele in \cite{TT03}, where the Fourier-analytic properties of the assignment $F \mapsto (a,b)$ were derived and this assignment was called the \emph{(discrete) nonlinear Fourier transform} or the \emph{discrete NFT} for short.
Its relationship to the orthogonal polynomials on the unit circle was also explained in \cite{TT03}; compare with the monographs by Simon \cite{S05a,S05b}. Let us only mention that
\[ \begin{bmatrix} e^{2\pi i n t} & 0 \\ 0 & e^{-2\pi i n t} \end{bmatrix}
\bigg( \prod_{n=1}^{N} \begin{bmatrix} A_n & B_n e^{2\pi i n t} \\ \overline{B_n} e^{-2\pi i n t} & A_n \end{bmatrix} \bigg)^{\textup{T}} \]
for positive integers $N$ are known as the \emph{Szeg\H{o} matrices} associated with the \emph{Verblunsky coefficients} $(-F_{n+1})_{n=0}^{\infty}$, and they are more commonly viewed as matrix functions of the complex variable $z=e^{2\pi i t}$.
The same model was also studied in a couple of recent papers by two of the present authors, \cite{O17} and \cite{R20}.
Discrete NFT can also be thought of as a discrete variant of the Dirac scattering transform; see the definition in \cite{MTT03a}.

{Discrete NFT is, just as its name already suggests, closely tied to the ordinary (linear) Fourier transform of double-sided complex sequences. Indeed, if the sequence $F$ is small (say, in the $\ell^1$-norm) then the defining formula \eqref{eq:defforNFT} easily gives that $b(t)$ is equal to $\widehat{F}(t) := \sum_{n\in\mathbb{Z}}F_n e^{2\pi i n t}$, up to a quadratic error in $F$. Consequently, we expect that the nonlinear transform $F\mapsto(a,b)$ still reflects some of the many useful properties of the linear Fourier transform $F\mapsto\widehat{F}$. This is in fact the case and the nonlinear analogues of well-known results in the Fourier analysis have received some attention over the last $20$ years.}

Let us briefly recall the basic {identities and estimates satisfied by} the discrete NFT.
A well-known identity, first formulated by Verblunsky \cite{V36}, in our notation reads:
\[ \int_\mathbb{T} \log|a(t)| dt = \sum_{n\in\mathbb{Z}} \log A_n = -\frac{1}{2} \sum_{n\in\mathbb{Z}} \log(1-|F_n|^2), \]
or more conveniently (multiplying both sides by $2$):
\[ \int_\mathbb{T} \log|a(t)|^2 dt = \sum_{n\in\mathbb{Z}} \log|A_n|^2, \]
i.e.
\begin{equation}\label{eq:p}
\big\|(\log|a(t)|^2)^{1/2}\big\|_{\textup{L}^2_t(\mathbb{T})} = \big\|(\log|A_n|^2)^{1/2}\big\|_{\ell^2_n(\mathbb{Z})}.
\end{equation}
It is {{thus}} regarded as the nonlinear analogue of Parseval's identity.

The \emph{nonlinear Hausdorff--Young inequality} was formulated in \cite{TT03} as:
\[ \big\|(\log|a(t)|)^{1/2}\big\|_{\textup{L}^q_t(\mathbb{T})} \leq \widetilde{C}_p \big\||\log(1-|F_n|^2)|^{1/2}\big\|_{\ell^p_n(\mathbb{Z})}. \]
Taking into account that $1-|F_n|^2=A_n^{-2}$ and denoting $C_p=2^{1/2}\widetilde{C}_p$, {this can be rewritten as}:
\begin{equation}\label{eq:hy}
\big\|(\log|a(t)|^2)^{1/2}\big\|_{\textup{L}^q_t(\mathbb{T})} \leq C_p \big\|(\log|A_n|^2)^{1/2}\big\|_{\ell^p_n(\mathbb{Z})}
\end{equation}
for $1\leq p<2$ and its conjugated exponent $2<q\leq\infty$, i.e., $1/p+1/q=1$.
Its proof was sketched in \cite{TT03} as an adaptation of the techniques by Christ and Kiselev \cite{CK01a,CK01b} to the discrete-parameter setting; also see \cite{O17} for the considerably stronger, nonlinear variational Hausdorff--Young inequality, which implies \eqref{eq:hy}.
However, the existing proofs give constants $C_p$ that blow up as $p\to2^-$.
It is unknown {whether} inequality \eqref{eq:hy} holds with a constant that is independent of $1\leq p\leq 2$, despite {the fact} that the endpoint case $p=2$ is controlled by equality \eqref{eq:p}, while the other endpoint case $p=1$ is easy (giving $C_p=1$) and it is widely discussed in \cite{TT03}.
This open problem, either in discrete or continuous parameter, seems to have been first posed by Muscalu, Tao, and Thiele \cite{MTT03a}, and it was subsequently popularized in papers by several authors \cite{K12,KOR17,O17}.

At the time of writing it is {likewise} open whether the same estimate holds with constant $1$:
\begin{equation}\label{eq:hywith1}
\big\|(\log|a(t)|^2)^{1/2}\big\|_{\textup{L}^q_t(\mathbb{T})} \leq  \big\|(\log|A_n|^2)^{1/2}\big\|_{\ell^p_n(\mathbb{Z})},
\end{equation}
{even though such a concrete inequality might theoretically be easier to disprove.}
Inequality \eqref{eq:hywith1}, if true, {would} be in perfect analogy with the classical sharp linear Hausdorff--Young inequality on $\mathbb{Z}$ (see \cite[Section XII.2, Theorem 2.3(ii)]{Z02}):
\begin{equation}\label{eq:linhy}
\big\|\widehat{G}\big\|_{\textup{L}^q(\mathbb{T})} \leq  \|G\|_{\ell^p(\mathbb{Z})}
\end{equation}
for a finitely supported sequence of complex numbers $G=(G_n)_{n\in\mathbb{Z}}$ and {its Fourier transform $\widehat{G}(t)=\sum_{n\in\mathbb{Z}} G_n e^{2\pi i n t}$, which is simply} the trigonometric polynomial having $G$ as its coefficients.

The purpose of this short note is to give some supporting evidence to the aforementioned conjecture on uniform boundedness of constants $C_p$ in \eqref{eq:hy} and {to} the conjectured sharp estimate \eqref{eq:hywith1}.
In the rest of the paper we always assume that $1<p<2$ and $2<q<\infty$ are given and that they are related by the H\"{o}lder scaling $1/p+1/q=1$.

Our first result discusses the constant $C_p$ in \eqref{eq:hy} for sequences $F$ with a small $\ell^1$-norm.

\begin{theorem}\label{thm:general}
If a sequence of coefficients $F$ satisfies $\|F\|_{\ell^1(\mathbb{Z})}\leq1/2$, then
\begin{equation}\label{eq:hyweaker}
\big\|(\log|a(t)|^2)^{1/2}\big\|_{\textup{L}^q_t(\mathbb{T})} \leq \big(1+3\|F\|_{\ell^1(\mathbb{Z})}\big) \,\big\|(\log|A_n|^2)^{1/2}\big\|_{\ell^p_n(\mathbb{Z})} .
\end{equation}
\end{theorem}

{The proof of Theorem~\ref{thm:general} will reveal that the numbers $1/2$ and $3$ in its formulation are somewhat arbitrary. One can increase the threshold for $\|F\|_{\ell^1(\mathbb{Z})}$ from $1/2$ to anything strictly smaller than $1$ at the cost of largely increasing the coefficient next to $\|F\|_{\ell^1(\mathbb{Z})}$ in \eqref{eq:hyweaker}. Similarly, the number $3$ in \eqref{eq:hyweaker} can be lowered to anything strictly larger than $1$ by being more restrictive on the magnitude of $\|F\|_{\ell^1(\mathbb{Z})}$.}

An immediate consequence of \eqref{eq:hyweaker} is a bound that is uniform in $p$:
\[ \big\|(\log|a(t)|^2)^{1/2}\big\|_{\textup{L}^q_t(\mathbb{T})} \leq \frac{5}{2} \big\|(\log|A_n|^2)^{1/2}\big\|_{\ell^p_n(\mathbb{Z})} , \]
as {long} as $\|F\|_{\ell^1(\mathbb{Z})}\leq1/2$.
However, the main point of the stronger estimate is that \eqref{eq:hyweaker} also gives an asymptotically sharp bound:
\[ \big\|(\log|a(t)|^2)^{1/2}\big\|_{\textup{L}^q_t(\mathbb{T})} \leq \big(1+O(\|F\|_{\ell^1(\mathbb{Z})})\big) \,\big\|(\log|A_n|^2)^{1/2}\big\|_{\ell^p_n(\mathbb{Z})}  \]
as  $\|F\|_{\ell^1(\mathbb{Z})}\to 0$.
Theorem~\ref{thm:general} is shown in Section~\ref{sec:general} below. Essentially, one only needs to estimate the error coming from linearizing \eqref{eq:hy}, but some care is needed if we want to deduce \eqref{eq:hyweaker} as it is formulated.

In our second result the emphasis is on the sharp constant $1$ in \eqref{eq:hywith1}, but only for sequences that satisfy an additional condition \eqref{eq:condition} below.

\begin{theorem}\label{thm:special}
There {exist} numbers $\alpha,\delta>0$ (depending on $1<p<2$) such that the following holds:
if a sequence of coefficients $F$ {is not identically $0$ and it} satisfies
\begin{equation}\label{eq:condition}
\|F\|_{\ell^1(\mathbb{Z})}\leq \delta \bigg(1-\frac{\|F\|_{\ell^\infty(\mathbb{Z})}}{\|F\|_{\ell^p(\mathbb{Z})}}\bigg)^\alpha ,
\end{equation}
then the corresponding $\textup{SU}(1,1)$-valued trigonometric product satisfies {the} sharp inequality \eqref{eq:hywith1}.
\end{theorem}

In Section~\ref{sec:nearextremizers} we recall the sharpened linear Hausdorff--Young inequality of Charalambides and Christ \cite{CK11}, which characterizes its near-extremizers, i.e., the sequences for which the (linear) Haus\-dorff--Young constant is close to $1$.
It will be a crucial ingredient in our proof of Theorem~\ref{thm:special} in Section~\ref{sec:special}.
The same discussion will also shed light on condition \eqref{eq:condition}, as we will recognize sequences {satisfying \eqref{eq:condition}} as being far from linear Hausdorff--Young extremizers (relative to their $\ell^1$-norm). In elementary terms, {such} sequences are sufficiently spread out over their support, as measured by the quantity $1-\|F\|_{\ell^\infty(\mathbb{Z})}/\|F\|_{\ell^p(\mathbb{Z})}$, relative to their size, as measured by $\|F\|_{\ell^1(\mathbb{Z})}$.

Typical instances of sequences $F$ that fail condition \eqref{eq:condition} are those that have only one nonzero term, as then the right hand side in \eqref{eq:condition} vanishes. The latter sequences {exhibit} the opposite behavior, as they are the exact extremizers of the linear Hausdorff--Young inequality. However, for these sequences, $t\mapsto a(t)$ is a constant function and \eqref{eq:hywith1} also holds, as it becomes a trivial equality. {We were not able to implement this dichotomy for the possible proof of \eqref{eq:hywith1} for all small sequences $F$, even though this could be a plausible strategy.}

Related work was done in the continuous {setting} by the present authors in \cite{KOR17}. That paper dealt with the nonlinear Hausdorff--Young inequality for ``small'' functions on the real line and the $\textup{SU}(1,1)$ (i.e.\@ the Dirac) scattering transform.
{Similarly as in \cite{KOR17}, the proof of Theorem~\ref{thm:special} actually ``beats'' the constant $1$ in \eqref{eq:hywith1} in a more quantitative way for sequences $F$ satisfying \eqref{eq:condition} and gives
\[ \big\|(\log|a(t)|^2)^{1/2}\big\|_{\textup{L}^q_t(\mathbb{T})} \leq \big(1-9\|F\|_{\ell^1(\mathbb{Z})}^2\big) \big\|(\log|A_n|^2)^{1/2}\big\|_{\ell^p_n(\mathbb{Z})}. \]
In particular the equality in \eqref{eq:hywith1} is never attained for such sequences.
We avoided formulating Theorem~\ref{thm:special} in that way, because the factor $1-9\|F\|_{\ell^1(\mathbb{Z})}^2$ is no longer optimal in any way.}

In the $\textup{SU}(1,1)$-valued setting there is no direct transference from the continuous model to the discrete one, so the discrete-parameter results often turn out more difficult to establish. For instance, the paper \cite{O17} is a discrete counterpart of an inequality by Oberlin, Seeger, Tao, Thiele, and Wright from \cite{OSTTW12}.

If one insists on studying general functions, {then the uniformity of constants $C_p$ in \eqref{eq:hy} is known to hold for a slightly different model of the nonlinear Fourier transform, where the exponentials $e^{2\pi i n t}$ in \eqref{eq:defforNFT} are replaced by characters of the so-called Cantor group model of the real line; see the paper \cite{K12} by one of the present authors. In that setting one gives up} the usual group structure and the topology of $\mathbb{R}$. {The Cantor group analogue was suggested in \cite{MTT03a}, where the authors discussed a toy model of a different problem for the nonlinear NFT.}

\section{Proof of Theorem~\ref{thm:general}}
\label{sec:general}
In this section we assume that $F$ is a nonzero sequence satisfying $\|F\|_{\ell^1(\mathbb{Z})}\leq1/2$.
Let us begin with a preliminary observation
\[ (\log|A_n|^2)^{1/2} = \big(-\log(1-|F_n|^2)\big)^{1/2} \geq |F_n|, \]
which in particular implies
\begin{equation}
\label{eq:auxseqest1}
\|F\|_{\ell^{p}(\mathbb{Z})} \leq \big\|(\log|A_n|^2)^{1/2}\big\|_{\ell^p_n(\mathbb{Z})}.
\end{equation}
The elementary inequality
\[ \prod_{n\in\mathbb{Z}} (1-x_n) \geq 1 - \sum_{n\in\mathbb{Z}} x_n \]
holds for arbitrary numbers $x_n\in[0,1]$, $n\in\mathbb{Z}$, {as  is} easily shown by {mathematical induction} on the number of nonzero summands and passing to the limit.
Using this inequality with $x_n=|F_n|^2$ we estimate
\[ \prod_{n\in\mathbb{Z}}A_n = \Big(\frac{1}{\prod_{n\in\mathbb{Z}}(1-|F_n|^2)}\Big)^{1/2}
\leq \Big(\frac{1}{1-\sum_{n\in\mathbb{Z}}|F_n|^2}\Big)^{1/2}
\leq \Big(\frac{1}{1-\|F\|_{\ell^{1}(\mathbb{Z})}^2}\Big)^{1/2}. \]
{Recalling the smallness condition on $F$ we}, in turn, obtain
\begin{equation}
\label{eq:auxseqest2}
\prod_{n\in\mathbb{Z}}A_n \leq 1 + \|F\|_{\ell^1(\mathbb{Z})}^2.
\end{equation}

{Introducing} the partial products
\[ \begin{bmatrix} a_N(t) & b_N(t) \\ \overline{b_N(t)} & \overline{a_N(t)} \end{bmatrix}
:= \prod_{n=-\infty}^{N} \begin{bmatrix} A_n & B_n e^{2\pi i n t} \\ \overline{B_n} e^{-2\pi i n t} & A_n \end{bmatrix}, \]
then from
\[ \begin{bmatrix} a_N(t) & b_N(t) \\ \overline{b_N(t)} & \overline{a_N(t)} \end{bmatrix}
= \begin{bmatrix} a_{N-1}(t) & b_{N-1}(t) \\ \overline{b_{N-1}(t)} & \overline{a_{N-1}(t)} \end{bmatrix}
\begin{bmatrix} A_N & B_N e^{2\pi i N t} \\ \overline{B_N} e^{-2\pi i N t} & A_N \end{bmatrix} \]
we immediately get {the} recurrence relations
\begin{align*}
a_N(t) & = a_{N-1}(t) A_N + b_{N-1}(t) \overline{B_N} e^{-2\pi i N t},\\
b_N(t) & = a_{N-1}(t) B_N e^{2\pi i N t} + b_{N-1}(t) A_N.
\end{align*}

The key idea is to introduce reduced quantities that will allow us to deduce a bootstrapping inequality for the $\textup{L}^q$-norms, which is \eqref{eq:bootstrap} below.
Denote
\[ \widetilde{a}_{N}(t) := \frac{a_N(t)}{\prod_{n=-\infty}^{N}A_n} - 1, \quad \widetilde{b}_{N}(t) := \frac{b_N(t)}{\prod_{n=-\infty}^{N}A_n}, \]
so that the recurrence relations become
\begin{align*}
\widetilde{a}_N(t) & = \widetilde{a}_{N-1}(t) + \widetilde{b}_{N-1}(t) \overline{F_N} e^{-2\pi i N t},\\
\widetilde{b}_N(t) & = \widetilde{b}_{N-1}(t) + F_N e^{2\pi i N t} + \widetilde{a}_{N-1}(t) F_N e^{2\pi i N t}
\end{align*}
and consequently give, by iteration,
\[ \widetilde{a}_{N}(t) = \sum_{n=-\infty}^{N-1} \widetilde{b}_{n}(t) \overline{F_{n+1}} e^{-2\pi i (n+1) t}, \quad
\widetilde{b}_{N}(t) = \sum_{n=-\infty}^{N} F_n e^{2\pi i n t} + \sum_{n=-\infty}^{N-1} \widetilde{a}_{n}(t) F_{n+1} e^{2\pi i (n+1) t}. \]
In particular,
\begin{equation}\label{eq:crucialiter}
\big|\widetilde{a}_{N}(t)\big| + \big|\widetilde{b}_{N}(t)\big|
\leq \sum_{n=-\infty}^{N-1} |F_{n+1}| \Big( \big|\widetilde{a}_{n}(t)\big| + \big|\widetilde{b}_{n}(t)\big| \Big)
+ \Big| \sum_{n=-\infty}^{N} F_n e^{2\pi i n t} \Big|,
\end{equation}
so that the linear Hausdorff--Young inequality applied to $\ldots, F_{N-1}, F_N, 0, 0, \ldots$ implies
\begin{equation}\label{eq:bootstrap}
\Big\| \big|\widetilde{a}_{N}\big| + \big|\widetilde{b}_{N}\big| \Big\|_{\textup{L}^q(\mathbb{T})}
\leq \sum_{n=-\infty}^{N-1} |F_{n+1}| \Big\| \big|\widetilde{a}_{n}\big| + \big|\widetilde{b}_{n}\big| \Big\|_{\textup{L}^q(\mathbb{T})}
+ \|F\|_{\ell^p(\mathbb{Z})}.
\end{equation}
Let $N_{\textup{min}}$ (resp.\@ $N_{\textup{max}}$) be the smallest (resp.\@ largest) integer $n$ such that $F_n\neq 0$.
By {mathematical induction} over $N\in\mathbb{Z}$, $N\geq N_{\textup{min}}-1$, estimate \eqref{eq:bootstrap} proves
\[ \Big\| \big|\widetilde{a}_{N}\big| + \big|\widetilde{b}_{N}\big| \Big\|_{\textup{L}^q(\mathbb{T})} \leq \frac{\|F\|_{\ell^p(\mathbb{Z})}}{1-\|F\|_{\ell^1(\mathbb{Z})}} \]
and, by taking $N=N_{\textup{max}}$, {gives}
\begin{equation}\label{eq:weakermain}
\big\|(\log|a|^2)^{1/2}\big\|_{\textup{L}^q(\mathbb{T})} \leq \|b\|_{\textup{L}^q(\mathbb{T})}
\leq \frac{\prod_{n\in\mathbb{Z}}A_n}{1-\|F\|_{\ell^1(\mathbb{Z})}} \|F\|_{\ell^p(\mathbb{Z})}.
\end{equation}
Finally, \eqref{eq:weakermain} in combination with \eqref{eq:auxseqest1} and \eqref{eq:auxseqest2} implies \eqref{eq:hyweaker}.

\section{Near-extremizers of the linear Hausdorff--Young inequality}
\label{sec:nearextremizers}
We will need the following result of Charalambides and Christ \cite[Theorem~1.3]{CK11}; also see \cite{Ch14}. It can be thought of as a simpler discrete variant of the analogous continuous-parameter result of Christ \cite{C14}. We only state it in dimension $d=1$.

\begin{theorem}[from \cite{CK11}]
\label{thm:CC}
For $1<p<2$ there exist constants $c,\gamma,\eta>0$ and a continuous nondecreasing function $\Lambda\colon\langle0,1]\to\langle0,1]$, all depending on $p$, such that
\[ \Lambda(t) \leq 1 - c(1-t)^\gamma \text{ for each } t\in[1-\eta,1] \]
and for all sequences $G$ {that are not identically $0$} we have
\[ { \big\| \widehat{G} \big\|_{\textup{L}^q(\mathbb{T})} } \leq \Lambda\Big(\frac{\|G\|_{\ell^\infty(\mathbb{Z})}}{\|G\|_{\ell^p(\mathbb{Z})}}\Big) \|G\|_{\ell^p(\mathbb{Z})} . \]
\end{theorem}

Exact extremizers $G$ of the linear Hausdorff--Young inequality \eqref{eq:linhy} are sequences with precisely one nonzero term.
Theorem~\ref{thm:CC} {states} that every ``near-extremizer'' $G$ of \eqref{eq:linhy} {must} be a sequence with the ratio $\|G\|_{\ell^\infty}/\|G\|_{\ell^p}$ close to $1$, which means that it is predominantly supported on a single point.
Theorem~\ref{thm:CC} is a highly nontrivial result {which} relies on several ideas from additive combinatorics; see \cite{CK11} or \cite{Ch14} for details.

\section{Proof of Theorem~\ref{thm:special}}
\label{sec:special}
Let $c,\gamma,\eta,\Lambda$ be as in Theorem~\ref{thm:CC}. We will take
\[ \alpha := \max\{1,\gamma\}, \quad \delta := \min\Bigg\{\frac{1}{6}, \frac{c\eta^\gamma}{3}, \bigg(3+\Big(\frac{3}{c}\Big)^{1/\gamma}\bigg)^{-\alpha} \Bigg\}, \]
and assume that $F$ satisfies condition \eqref{eq:condition}. In particular,
\[ \|F\|_{\ell^1(\mathbb{Z})}\leq\delta \leq \frac{1}{6} \]
and
\begin{equation}\label{eq:simplealgebra}
3\|F\|_{\ell^1(\mathbb{Z})} + \bigg(\frac{3\|F\|_{\ell^1(\mathbb{Z})}}{c}\bigg)^{1/\gamma}
\leq 3\|F\|_{\ell^1(\mathbb{Z})}^{1/\alpha} + \bigg(\frac{3}{c}\bigg)^{1/\gamma} \|F\|_{\ell^1(\mathbb{Z})}^{1/\alpha}
\leq \bigg(\frac{\|F\|_{\ell^1(\mathbb{Z})}}{\delta}\bigg)^{1/\alpha}.
\end{equation}

Note that \eqref{eq:weakermain} from Section~\ref{sec:general} is insufficient for the proof of sharp inequality \eqref{eq:hywith1}, as the right hand side of \eqref{eq:weakermain} can be larger than the right hand side of \eqref{eq:hywith1}.
(Simply take any $F$ with two nonzero terms $F_1,F_2$ such that $0<|F_1|=|F_2|<1/2$.)
An improvement will be obtained if we postpone the application of the linear Hausdorff--Young inequality and from \eqref{eq:crucialiter} only conclude
\[ \Big\| \big|\widetilde{a}_{N}\big| + \big|\widetilde{b}_{N}\big| \Big\|_{\textup{L}^q(\mathbb{T})}
\leq \sum_{n=-\infty}^{N-1} |F_{n+1}| \Big\| \big|\widetilde{a}_{n}\big| + \big|\widetilde{b}_{n}\big| \Big\|_{\textup{L}^q(\mathbb{T})}
+ \sup_{N\in\mathbb{Z}}\Big\| \sum_{n=-\infty}^{N} F_n e^{2\pi i n t} \Big\|_{\textup{L}^q_t(\mathbb{T})}, \]
which, using induction in the same way as before, leads to
\begin{equation}\label{eq:strongermain}
\big\|(\log|a|^2)^{1/2}\big\|_{\textup{L}^q(\mathbb{T})} \leq \big(1+3\|F\|_{\ell^1(\mathbb{Z})}\big)\, \sup_{N\in\mathbb{Z}}\Big\| \sum_{n=-\infty}^{N} F_n e^{2\pi i n t} \Big\|_{\textup{L}^q_t(\mathbb{T})} .
\end{equation}
We claim that condition \eqref{eq:condition} implies
\begin{equation}\label{eq:sharperlinhy}
\Big\| \sum_{n=-\infty}^{N} F_n e^{2\pi i n t} \Big\|_{\textup{L}^q_t(\mathbb{T})} \leq \big(1-3\|F\|_{\ell^1(\mathbb{Z})}\big) \,\|F\|_{\ell^p(\mathbb{Z})}
\end{equation}
for each $N\in\mathbb{Z}$, which then combines with \eqref{eq:auxseqest1} and \eqref{eq:strongermain} to establish \eqref{eq:hywith1}.
Thus, it remains to prove \eqref{eq:sharperlinhy}.

For a given $N\in\mathbb{Z}$ denote the truncated sequence:
\[ \widetilde{F} := (\ldots, F_n, \ldots, F_{N-1}, F_N, 0, 0, \ldots ). \]
We distinguish two possible cases in terms of $N$.

$(1^\circ)$
If $N$ is (sufficiently small) such that
\[ \|\widetilde{F}\|_{\ell^p(\mathbb{Z})} \leq \big(1-3\|F\|_{\ell^1(\mathbb{Z})}\big) \,\|F\|_{\ell^p(\mathbb{Z})} , \]
then we apply the ordinary linear Hausdorff--Young inequality \eqref{eq:linhy} to $\widetilde{F}$:
\[ \Big\| \sum_{n=-\infty}^{N} F_n e^{2\pi i n t} \Big\|_{\textup{L}^q_t(\mathbb{T})} \leq \|\widetilde{F}\|_{\ell^p(\mathbb{Z})} \]
and it turns precisely into \eqref{eq:sharperlinhy}.

$(2^\circ)$
If $N$ is (sufficiently large) such that
\[ \|\widetilde{F}\|_{\ell^p(\mathbb{Z})} > \big(1-3\|F\|_{\ell^1(\mathbb{Z})}\big) \,\|F\|_{\ell^p(\mathbb{Z})} , \]
then condition \eqref{eq:condition} followed by \eqref{eq:simplealgebra} gives
\[ \frac{\|\widetilde{F}\|_{\ell^\infty(\mathbb{Z})}}{\|\widetilde{F}\|_{\ell^p(\mathbb{Z})}} \leq \frac{1}{1-3\|F\|_{\ell^1(\mathbb{Z})}} \frac{\|F\|_{\ell^\infty(\mathbb{Z})}}{\|F\|_{\ell^p(\mathbb{Z})}}
\leq \frac{1-(\|F\|_{\ell^1(\mathbb{Z})}/\delta)^{1/\alpha}}{1-3\|F\|_{\ell^1(\mathbb{Z})}} \leq 1 - \bigg(\frac{3\|F\|_{\ell^1(\mathbb{Z})}}{c}\bigg)^{1/\gamma}, \]
which implies
\[ \Lambda\Big(\frac{\|\widetilde{F}\|_{\ell^\infty(\mathbb{Z})}}{\|\widetilde{F}\|_{\ell^p(\mathbb{Z})}}\Big) \leq \max\{1-3\|F\|_{\ell^1(\mathbb{Z})}, 1-c\eta^{\gamma}\} = 1-3\|F\|_{\ell^1(\mathbb{Z})}. \]
This time \eqref{eq:sharperlinhy} is guaranteed by Theorem~\ref{thm:CC} applied to $\widetilde{F}$. The proof is now complete.


\section*{Acknowledgements}
VK and JR were supported in part by the \emph{Croatian Science Foundation} project UIP-2017-05-4129 (MUNHANAP).
DOS was supported by the \emph{Engineering and Physical Sciences Research Council} New Investigator Award ``Sharp Fourier Restriction Theory'', grant no.\@ EP/T001364/1.

We are grateful to the careful referee for numerous comments and suggestions.


\bibliography{discrete_nonlinear_hy}{}

\begin{thebibliography}{10}

\bibitem{CK11}
Marcos {Charalambides} and Michael {Christ}.
\newblock {Near-extremizers of Young's inequality for discrete groups}.
\newblock Preprint available at arXiv:1112.3716, 2011.

\bibitem{Ch14}
Marcos {Charalampidis}.
\newblock {\em Extremal Problems in Analysis}.
\newblock PhD thesis, University of California, Berkeley, 2014.

\bibitem{C14}
Michael {Christ}.
\newblock {A sharpened Hausdorff-Young inequality}.
\newblock Preprint available at arXiv:1406.1210, 2014.

\bibitem{CK01a}
Michael {Christ} and Alexander {Kiselev}.
\newblock {Maximal functions associated to filtrations}.
\newblock {\em {J. Funct. Anal.}}, 179(2):409--425, 2001.

\bibitem{CK01b}
Michael {Christ} and Alexander {Kiselev}.
\newblock {WKB asymptotic behavior of almost all generalized eigenfunctions for
  one-dimensional Schr\"odinger operators with slowly decaying potentials}.
\newblock {\em {J. Funct. Anal.}}, 179(2):426--447, 2001.

\bibitem{K12}
Vjekoslav {Kova\v{c}}.
\newblock {Uniform constants in Hausdorff-Young inequalities for the Cantor
  group model of the scattering transform}.
\newblock {\em {Proc. Amer. Math. Soc.}}, 140(3):915--926, 2012.

\bibitem{KOR17}
Vjekoslav {Kova\v{c}}, Diogo {Oliveira e Silva}, and Jelena {Rup\v{c}i\'c}.
\newblock {A sharp nonlinear Hausdorff-Young inequality for small potentials}.
\newblock {\em {Proc. Amer. Math. Soc.}}, 147(1):239--253, 2019.

\bibitem{MTT03a}
Camil {Muscalu}, Terence {Tao}, and Christoph {Thiele}.
\newblock {A Carleson type theorem for a Cantor group model of the scattering
  transform}.
\newblock {\em {Nonlinearity}}, 16(1):219--246, 2003.

\bibitem{OSTTW12}
Richard {Oberlin}, Andreas {Seeger}, Terence {Tao}, Christoph {Thiele}, and
  James {Wright}.
\newblock {A variation norm Carleson theorem}.
\newblock {\em {J. Eur. Math. Soc.}}, 14(2):421--464, 2012.

\bibitem{O17}
Diogo {Oliveira e Silva}.
\newblock {A variational nonlinear Hausdorff-Young inequality in the discrete
  setting}.
\newblock {\em {Math. Res. Lett.}}, 25(6):1993--2015, 2018.

\bibitem{R20}
Jelena {Rup\v{c}i\'c}.
\newblock {Convergence of lacunary SU(1,1)-valued trigonometric products}.
\newblock {\em {Commun. Pure Appl. Anal.}}, 19(3):1275--1289, 2020.

\bibitem{S05a}
Barry {Simon}.
\newblock {\em {Orthogonal polynomials on the unit circle. Part 1: Classical
  theory}}, volume~54.
\newblock AMS, Providence, 2005.

\bibitem{S05b}
Barry {Simon}.
\newblock {\em {Orthogonal polynomials on the unit circle. Part 2: Spectral
  theory}}, volume~54.
\newblock AMS, Providence, 2005.

\bibitem{TT03}
Terence {Tao} and Christoph {Thiele}.
\newblock {Nonlinear Fourier Analysis}.
\newblock Unpublished lecture notes available at arXiv:1201.5129, 2003.

\bibitem{V36}
Samuel {Verblunsky}.
\newblock {On positive harmonic functions. II}.
\newblock {\em {Proc. London Math. Soc. (2)}}, 40:290--320, 1935.

\bibitem{Z02}
Antoni {Zygmund}.
\newblock {\em Trigonometric series. Vols. 1 \& 2, Third Ed.}
\newblock Cambridge University Press, 2002.

\end{thebibliography}
\bibliographystyle{plain}

\end{document}